\documentclass[english,russian]{amsart}
\usepackage{amssymb,amsmath}

\newtheorem{thm}{Theorem}[section]

\newtheorem{lm}[thm]{Lemma}

\newcommand{\sudda}[1]{}

\begin{document}

\title{Symmetry groups of pfaffians of symmetric matrices}

\author{A.S. Dzhumadil'daev}

\address{Kazakh-British Technical  University, Tole bi 59,
Almaty, 050000,Kazakhstan}
\email{dzhuma@hotmail.com}

\begin{abstract} We prove that symmetry group of the pfaffian polynomial of a symmetric matrix is a dihedral group. We calculate pfaffians of symmetric matrices  with components $(x_i-x_j)^2$ and $\cos(x_i-x_j)$
for $i<j.$\end{abstract}




\maketitle

\section{Introduction}

In first glance the title of the paper is not correct. Pfaffians usually are connected with determinants of skew-symmetric matrices. If $a_{i,j}=-a_{j,i},$ for any $1\le i,j\le 2n,$ then determinant of skew-symmetric matrix $A=(a_{i,j})$
is complete square and square root of determinant is pfaffian,
$$det\,A=(pf_{2n} A)^2.$$
In fact, pfaffian is defined not for whole  matrix $A.$ To construct pfaffians it is enough to know upper triangular part of $A.$  

Connection between determinants of  skew-symmetric matrices and pfaffians was firstly noted in \cite{Cayley}.  Details of pfaffian constructions see also \cite{Bourbaki} and \cite{Vein}.

Let $S_{2n}$ be set of permutations on the set $[2n]=\{1,2,\ldots,2n\}$ and $S_{2n,pf}$  its subset of permutations called {\it Pfaff permutations},
$$S_{2n,pf}=\{ \sigma=(i_1,j_1,\ldots,i_n, j_n)\in S_{2n} | i_1<i_2<\cdots <i_n, i_s<j_s, 1\le s\le n\}.$$ 
For any $\sigma\in S_{2n,pf}$ we define {\it Pffaf aggregates} $a_{\sigma}$ by 
$$a_{\sigma}=a_{i_1,j_1}\cdots a_{i_n,j_n}.$$
We see that pfaffian aggregates  are defined for any triangular array $\bar A=(a_{i,j})_{1\le i<j\le 2n}.$ 
Then pffaffian of order $2n$ is a polynomial defined as an alternating sum of pfaff aggregates
$$pf_n=\sum_{\sigma\in S_{2n,pf}} sign\;\sigma\; a_{\sigma}.$$

Now suppose  that $\{a_{i,j}, 1\le i,j\le 2n\}$  are $n^2$ generators and endow space polynomials 
$K[a_{i,j} | 1\le i,j\le n]$ by structure of $S_{2n}$-module by the following action on generators 
$$\sigma a_{i,j}=a_{\sigma^{-1}(i),\sigma^{-1}(j)}.$$
In particular, if $A=(a_{i,j})$ with skew-symmetric set of generators, $a_{i,j}=-a_{j,i}$ then this action induces 
structure of $S_{2n}$-module structure on the space of polynomials with ${n\choose 2}$ generators
$K[a_{i,j} | 1\le i<j\le n].$ Similarly, we obtain one more  structure of $S_{2n}$-module on this space if generators are symmetric, $a_{i,j}=a_{j,i}.$ In both cases appear natural questions about invariants under these actions of permutation groups. In particular, we can ask about symmetry  and skew-symmetry groups
of given polynomial   $f\in K[a_{i,j}],$
$$Sym\; f=\{\sigma\in S_n | \sigma \; f=f\},$$
$$SSym\; f=\{\sigma\in S_n | \sigma \; f=sign \,\sigma \;f\}.$$

For example,  determinant polynomial $det\,A $ for $A=(a_{ij})_{1\le i,j\le n}$  is a polynomial of degree $n$ and its symmetry group is isomorphic to $S_n.$  Another example: if matrix $A$ is skew-symmetric, then pfaff polynomiaml$pf_{2n}=pf_{2n}A$ is a polynomial of degree $n$ and 
$$SSym \; pf_{2n}\cong S_{2n}.$$

Let characteristic of main field is $p \ne 2.$ The aim of our paper is to prove the following result.

\begin{thm} \label{main}  If generators $a_{i,j}$ are symmetric, $a_{i,j}=a_{j,i},$ then symmetry  group of pffaffian polynomial $pf_{2n}=pf_{2n} \bar A$ is isomorphic to a dihedral group
$$Sym \; pf_{2n} \cong D_{2n}.$$
\end{thm}

 In the proof of this fact we use the following results.

\begin{thm} \label{translation} 
Let $\psi$ be symmetric function, 
$$\psi(x,y)=\psi(y,x),$$
such that
\begin{equation}\label{const}
\psi(x+z,y+z)=\psi(x,y),
\end{equation}
for any $x,y,z.$ 
Let $pf_{2n}(x_1,\ldots,x_{2n})$ be pfaffian polynomial $pf_{2n} \bar A,$ where $a_{i,j}=\psi(x_i,x_j).$
Then for any $1\le s\le 2n,$
$$ pf_{2n} (x_1,\ldots,x_{s-1}, x_s,x_s, x_{s+2},\ldots,x_{2n})=
 c \; pf_{2n-2} (x_1,\ldots,x_{s-1},x_{s+2},\ldots,x_{2n}), $$
where 
$$c=\psi(0,0).$$
\end{thm}
Here we set $x_{2n+1}=x_1.$

\begin{thm} \label{2222} Let $\bar A=(a_{i,j})_{1\le i<j \le 2n}$ be triangular array with components  $a_{i,j}=(x_i-x_j)^2$ for $1\le i<j\le 2n.$ Then 
$$pf_{2n}\;\bar A=-(-2)^{n-1}\; (x_1-x_2)(x_2-x_3)\cdots (x_{2n-1}-x_{2n})(x_{2n}-x_1).$$
\end{thm}

We give one more result on symmetric pfaffians.    

\begin{thm}\label{cos} If $a_{i,j}=cos(x_i-x_j),$ $1\le i<j\le 2n,$ then 
$$pf_{2n} A=\cos (x_1-x_2+x_3-x_4+\cdots +x_{2n-1}-x_{2n}).$$
\end{thm}

\section{ Proof of Theorem \ref{translation}}

Let $\bar A=(a_{i,j})_{1\le i<j\le 2n}$ be triangular array and 
$$hk_s \bar A=(a_{1,s}, a_{2,s},\ldots, a_{s-1,s},0,a_{s,s+1},\ldots,a_{s,2n})$$
be $s$-th hook of triangle $\bar A.$ 
If generators are skew-symmetric, then 
we can present it as $s$-th row of  
 skew-symmetric matrix $A,$ constructed by triangular $\bar A,$ 
$$hk_s A=(-a_{s,1}, -a_{s,2},\ldots, -a_{s,s-1},0,a_{s,s+1},\ldots,a_{s,2n})$$
Then decomposition formula for pfaffians by $s$-th hook can be written as follows
$$pf_{2n} A=\sum_{1\le j\le 2n, j \ne s} (-1)^{s+j+1+H(s-j)} a_{s,j} pf_{2n-2} A_{\hat{s}.\hat{j}},$$
where $A_{\hat{i},\hat{j}}$ denotes the matrix $A$ with both the $i$-th and $j$-th rows 
and columns removed.  
Here $\theta(i-j)$ is Heaviside step function,
$$H(s)=1, \mbox{ if $s>0$ and $0,$ if $s\le 0.$}$$

If generators $a_{i,j}$ are symmetric, $a_{i,j}=a_{j,i},$ then decomposition formula can be written without Heaviside function,
\begin{equation}\label{decomp}
pf_{2n} \bar A=\sum_{1\le j\le 2n, j \ne s} (-1)^{s+j+1} a_{s,j} pf_{2n-2} \bar A_{\hat{s},\hat{j}}.
\end{equation}

{\bf Proof of Theorem \ref{translation}.}  Take  $x=y=0$ in  (\ref{const}) . We obtain
\begin{equation}\label{const1}
\psi(z,z)=\psi(0,0)=c.
\end{equation}
By   (\ref{decomp})  we can  decompose pfaffian by $s$-th hook
\begin{equation}\label{1111}
  pf_{2n} \bar A=
  \end{equation}
  $$
 \sum_{i=1}^{s-1} (-1)^{s+i+1} a_{s,i} \; pf \; \{1,2,\ldots,\hat{i}\ldots,\hat{s}\ldots,2n\}+$$ $$
 a_{s,s+1} \; pf \; \{1,2,\ldots,\hat{s},\widehat{s+1},\ldots,2n\}+
 $$
 $$
 \sum_{i=s+2}^{2n}(-1)^{s+i+1}a_{s,i} \;pf \; \{ 1,2,\ldots, \hat{s}\ldots,\hat{i},\ldots,2n\}.
 $$
 Here notation $\hat{i}$ means that corresponding index $i$ is omitted. 
 
 Similarly decomposition of pfaffian by $s+1$-th hook gives us the following relation
  \begin{equation}\label{1112}
  pf_{2n} \bar A=
   \end{equation} 
   $$
 \sum_{i=1}^{s-1} (-1)^{s+i} a_{s+1,i} \;pf \; \{1,2,\ldots,\hat{i}\ldots,\widehat{s+1},\ldots,2n\}+$$
 $$
 a_{s+1,s} \;pf \; \{1,2,\ldots,\hat{s},\widehat{s+1},\ldots,2n\}+$$
 $$
 \sum_{i=s+2}^{2n}(-1)^{s+i}a_{s+1,i}\; pf \; \{ 1,2,\ldots, \widehat{s+1},\ldots,\hat{i},\ldots,2n\}.
 $$
 Take here $x_{s+1}=x_s.$ Then $a_{s+1,i}=a_{s,i},$ for any $i\ne s$ and by (\ref{const1}) 
 $$a_{s+1,s}=a_{s,s+1}=\psi(x_s,x_s)=\psi(0,0)=c.$$ 
 Therefore by  (\ref{1112})
  $$pf_{2n}(x_1,\ldots,x_{s-1},x_s,x_s,x_{s+2},\ldots,x_{2n})=$$
 $$\sum_{i=1}^{s-1} (-1)^{s+i} a_{s+1,i}\; pf\;  \{1,2,\ldots,\hat{i}\ldots,\widehat{s+1}\ldots,2n\}+$$
 $$\psi(0,0) \;pf\; \{ 1,2,\ldots, \hat{s},\widehat{s+1},\ldots,\ldots,2n\}+$$
 $$ \sum_{i=s+2}^{2n}(-1)^{s+i}a_{s,i} \;pf\; \{ 1,2,\ldots, \widehat{s+1},\ldots,\hat{i},\ldots,2n\}=$$
 $$\sum_{i=1}^{s-1} (-1)^{s+i} a_{s,i}\; pf\; \{1,2,\ldots,\hat{i}\ldots,\widehat{s+1},\ldots,2n\}+
 $$
  $$\psi(0,0) \;pf\; \{ 1,2,\ldots, \hat{s},\widehat{s+1},\ldots,\ldots,2n\}+$$
 $$
 \sum_{i=s+2}^{2n}(-1)^{s+i}a_{s,i} \;pf\; \{ 1,2,\ldots, \widehat{s+1},\ldots,\hat{i},\ldots,2n\}.$$
 Hence by  (\ref{1111})
  $$pf_{2n}(x_1,\ldots,x_{s-1},x_s,x_s,x_{s+2},\ldots,x_{2n})=$$
  
   $$\psi(0,0) \;pf\; \{ 1,2,\ldots, \hat{s},\widehat{s+1},\ldots,\ldots,2n\}+$$
 $$\sum_{i=1}^{s-1} (-1)^{s+i} a_{s,i}\; pf \;\{1,2,\ldots,\hat{i}\ldots,\hat{s}\ldots,2n\}+
 \sum_{i=s+2}^{2n}(-1)^{s+i}a_{s,i}\; pf\; \{ 1,2,\ldots, \hat{s}\ldots,\hat{i},\ldots,2n\}=$$
 
  $$\psi(0,0) \;pf\; \{ 1,2,\ldots, \hat{s},\widehat{s+1},\ldots,\ldots,2n\}+$$
 $$- \sum_{i=1}^{s-1} (-1)^{s+i+1} a_{s,i}\; pf\; \{1,2,\ldots,\hat{i}\ldots,\hat{s}\ldots,2n\}-
 \sum_{i=s+2}^{2n}(-1)^{s+i+1}a_{s,i} \; pf \; \{ 1,2,\ldots, \hat{s}\ldots,\hat{i},\ldots,2n\}=$$
 
  $$2\psi(0,0) \;pf\; \{ 1,2,\ldots, \hat{s},\widehat{s+1},\ldots,\ldots,2n\}- pf_{2n}(x_1,\ldots,x_{s-1},x_s,x_s,x_{s+2},\ldots,x_{2n}).$$
 So, 
  $$pf_{2n}(x_1,\ldots,x_{s-1},x_s, x_s,x_{s+2},\ldots,x_{2n})=
   \psi(0,0) \;pf\; \{ 1,2,\ldots, \hat{s},\widehat{s+1},\ldots,\ldots,2n\},$$
 if characteristic of main field is not $2.$ 
 $\square$

\section{Pfaffian of  $(x_i-x_j)^2$}

If $A=(a_{i,j})$ is skew-symmetric, then
$$det\,A=(pf_{2n} \bar A)^2.$$
If matrix $A$ is not skew-symmetric, say if $A$ is symmetric, then determinant polynomial $det\,A$ and pfaffian polynomial $pf_{2n}(\bar A)$ have no such connection. For example, if $A_{n}=((x_i-x_j)^2)_{1\le i,j\le n},$ then
$$det\,A_n=\left\{\begin{array}{ll}
-(x_1-x_2)^2,&\mbox{ if $n=2$,}\\
2((x_2-x_2)(x_2-x_3)(x_3-x_1))^2,&\mbox{ if $n=3$ and,}\\
0,&\mbox{ otherwise}\end{array}\right.$$
in that time, by Theorem \ref{2222},  pfaffians are non-trivial for any even $n.$

{\bf Proof of Theorem  \ref{2222}.}  Let $\psi(x,y)=(x-y)^2$ and $a_{i,j}=\psi(x_i,x_j)=(x_i-x_j)^2.$  Note that pfaffian $pf_{2n} A$ is a polynomial on variables $x_1,\ldots,x_{2n}.$ 
 Let us construct polynomial $g_{2n}=g_{2n}(x_1,x_2,\ldots,x_{2n})$ by 
$$g_{2n}=(x_1-x_2)(x_2-x_3)\cdots (x_{2n-1}-x_{2n})(x_{2n}-x_1).$$
We have to prove that 
$$pf_{2n}=-(-2)^{n-1}g_{2n}.$$

Since $\psi(0,0)=0,$  by Theorem \ref{translation}  the polynomial $pf_{2n}(x_1,\ldots,x_{s},x_{s+1},\ldots,x_{2n})$ is divisable by $x_{s}-~x_{s+1}$ for any $1\le s\le 2n.$  Here we set $x_{2n+1}=x_1.$ Note that degree of the polynomial
$g_{2n}(x_1,\ldots,x_{2n})$ is $2n$ and degree of $pf_{2n}((x_i-x_j)^2) $ is also $2n.$  
 Therefore,
 $$pf_{2n}(x_1,x_2,\ldots,x_{2n})=c\;g_{2n}(x_1,x_2,\ldots,x_{2n}),$$
 for some constant $c.$
Take $x_i=i.$   It is easy to see that 
  $$g_{2n} (1,2,\ldots, 2n)=
(1-2)(2-3)\cdots (2n-1-2n)(2n-1)=  
-(2n-1).$$
  
It remains to prove that 
\begin{equation}\label{1113}  pf_{2n}(1,2,\ldots, {2n})=(-2)^{n-1}(2n-1).
\end{equation}
to obtain $c=-(-2)^{n-1}.$

By induction on $n$ we will prove that 
 $$pf_{2n}(x_1,x_2,\ldots,x_{2n})=-(-2)^{n-1}g_{2n}(x_1,x_2,\ldots,x_{2n}).$$
 
 For $n=1$ our statement is evident:
   $$pf_2 \bar A=a_{1,2}=-(x_1-x_2)(x_2-x_1).$$
   
 Suppose that our statement is true for $n-1,$
$$pf_{2n-2} \bar A=
-(-2)^{n-2}(x_1-x_2)(x_2-x_3)\cdots (x_{2n-3}-x_{2n-2})(x_{2n-2}-x_1).$$
Let us  prove it for $n.$
   
We decompose pfaffian by the first row. By (\ref{decomp}) 
$$pf_{2n} \bar A=\sum_{i=2}^{2n} (-1)^{i} a_{1,i} \;pf_{2n-2} A_{\hat{1},\hat{i}}.$$
We see that 
$$pf_{2n}=R_1+R_2+R_3,$$
where 
$$R_1=a_{1,2}\; pf_{2n} \bar A_{\hat{1}\hat{2}},$$

$$R_2=\sum_{i=3}^{2n-1} (-1)^{i} a_{1,i} \; pf_{2n-2} \bar A_{\hat{1},\hat{i}},$$

$$R_3=a_{1,2n}\; pf_{2n-2}\bar A_{\hat{1}, \hat{2n}}.$$

By inductive suggestion
$$R_1=-(-2)^{n-2}(x_1-x_2)^2 (x_3-x_4)\cdots (x_{2n-1}-x_{2n})(x_{2n}-x_3).$$
Hence
$$R_1|_{x_i\rightarrow i}=
(1-2)^2(3-4)\cdots(2n-1-2n)(2n-3)=$$ 

$$(-1)^{2n-3}(2n-3)=(-2)^{n-2}(2n-3),$$

$$R_3|_{x_i\rightarrow i}=-(-2)^{n-2}(1-2n)^2(2-3)(3-4)\cdots (2n-2-2n+1)(2n-1-2)=$$

$$-(-2)^{n-2}(2n-1)^2(-1)^{2n-3}(2n-3)=
(-2)^{n-2}(2n-3)(2n-1)^2.$$
Further, if $2<i<2n,$ then 
$$(-1)^{i} a_{1,i} \; pf_{2n} \bar A_{\hat{1},\hat{i}} |_{x_j\rightarrow j}= $$

$$(-1)^i (-(-2)^{n-2}) (x_1-x_i)^2 (x_2-x_3)(x_3-x_4)\cdots (x_{i-1}-x_{i+1})(x_{i+1}-x_{i+2})\times \cdots $$
$$\times  (x_{2n-1}-x_{2n})(x_{2n}-x_2)|_{x_j\rightarrow j}=$$ 

$$
(-1)^i (- (-2)^{n-2} )(i-1)^2 (-2)(2n-2)=(-1)^i (i-1)^2 (-2)^{n-2}4(n-1).$$
Hence
$$R_2|_{x_i\rightarrow i}=-(-2)^{n-2}\sum_{i=3}^{2n-1} -(-1)^i (i-1)^2 4(n-1)= (-2)^{n-2} 4(n-1)(2n^2-3n+2).$$

So, we see that (\ref{1113}) is true for $n,$
$$f_{2n}(1,2,\ldots,2n)=R_1+R_2+R_3=$$

$$(-2)^{n-2}[(2n-3)-4(n-1)(2n^2-3n+2)+(2n-3)(2n-1)^2]=$$ 

$$
-(-2)^{n-2}2(2n-1)=(-2)^{n-1} (2n-1).$$

Theorem \ref{2222} is proved.  
$\square$

\section{Symmetry group of the polynomial $g_{2n}$ }

\begin{lm} \label{sym g}Let 
$$g_{2n}(x_1,\ldots, x_{2n})=(x_1-x_2)(x_2-x_3) \cdots (x_{2n-1}-x_{2n})(x_{2n}-x_1)$$
as in the proof of Theorem \ref{2222} and 
$$Sym(g_{2n})=\{\sigma\in S_{2n} | g_{2n}(x_{\sigma(1)},\ldots,x_{\sigma(2n)})=g_{2n}(x_1,\ldots x_{2n})\}$$
symmetric group of the polynomial $g_{2n}.$ 
Then 
$$Sym(g_{2n})\cong D_{2n}.$$
\end{lm}

{\bf Proof.}
First we  check that any dihedral permutation $\sigma\in D_{2n}$ is symmetry of the polynomial $g_{2n}.$ 

Let us take realization of dihedral group as a symmetry group of regular $n$-gon whose vertices are clockwise 
labelled by $1,2,\ldots,2n.$
Elements of dihedral group might have:
\begin{itemize}
\item[{\bf I.}] one up-run: $\sigma=1\;2\;\cdots (2n).$
\item[{\bf II}.] one down-run: $\sigma=(2n)\;(2n-1)\;\cdots \;1.$
\item[{\bf III.}] two up-runs 
$$\sigma(1)=s<\sigma(2)=s+1<\cdots <\sigma(2n-s+1)=2n, \;\; \sigma(2n-s+2)=1<\cdots <\sigma(2n)=s-1,$$
for some $1<s\le 2n,$
\item[{\bf IV}.] or two down-runs
$$\sigma(1)=s>\sigma(2)=s-1>\cdots >\sigma(s)=1,\;\; \sigma(s+1)=2n>\cdots >\sigma(2n)=s+1.$$
for some $1\le s<2n.$
\end{itemize}

In cases {\bf I} and {\bf II} our statement is evident.

In case {\bf III} we have 
$$g_{2n}(x_{\sigma(1)},\ldots,x_{\sigma(2n)})=$$
$$(x_s-x_{s+1})(x_{s+1}-x_{s+2})\cdots (x_{2n-1}-x_{2n})(x_{2n}-x_1)\; (x_1-x_2)\cdots(x_{s-2}-x_{s-1})(x_{s-1}-x_s)=$$
$$(x_1-x_2)\cdots (x_{2n-1}-x_{2n})(x_{2n}-x_1)=g_{2n}(x_1,\ldots,x_{2n}).$$

In  case {\bf IV} 
$$g_{2n}(x_{\sigma(1)},\ldots,x_{\sigma(2n)})=$$
$$(x_s-x_{s-1})(x_{s-1}-x_{s-2})\cdots (x_{2}-x_1)(x_1-x_{2n})\;(x_{2n}-x_{2n-1})\cdots(x_{s+1}-x_{s})=$$

$$(-1)^{s}(x_{s-1}-x_{s})(x_{s-2}-x_{s-1})\cdots (x_{1}-x_2)(x_{2n}-x_{1})\;
(-1)^{2n-s}(x_{2n-1}-x_{2n})\cdots(x_{s}-x_{s+1})=$$

$$(x_1-x_2)\cdots (x_{2n-1}-x_{2n})=g_{2n}(x_1,\ldots,x_{2n}).$$
So,
$$D_{2n}\subseteq Sym(g_{2n}).$$

Now we will prove that any $\sigma\in Sym(g_{2n})$ is a dihedral permutation. 

Let $M_{2n}=\{1,2,\ldots,2n\}.$ 
For $i,j\in M_{2n}$ say that they are connected, if $|i-j|=1$ or $|i-j|=2n-1.$
So, if  $i<j<2n,$ then $i,j$  are connected iff $j=i+1.$
If $j=2n,$ and $i,j$ are connected, then $i=2n-1$ or $i=1.$ 
It is clear that this relation is symmetric: $i,j$ are connected iff $j,i$ are connected. So,  $i,j\in M_{2n}$  are connected, if $|i-j|=1$ or
 $(i,j)=(1,2n) $ or $(i,j)=(2n,1).$
 
Note that he polynomial $g_{2n}(x_1,\ldots,x_{2n})$ is a product of polynomials $x_i-x_j,$ $i<j,$ where  $i$ and $j$ are connected. 
Therefore, any symmetry $\sigma\in Sym(g_{2n})$ has the following property: if $i$ and $j$ are connected, then $\sigma(i)$ and $\sigma(j)$ are also connected.

 Let $\sigma\in Sym(g_{2n})$ and $\sigma(1)=i_1.$ Might happen the following possibilities.
 
{\bf Case A.} Suppose that $\sigma(1)=i_1<\sigma(2).$ Take $k>1,$ such that $\sigma(k-1)<\sigma(k)$ and $\sigma(k+1)<\sigma(k).$ 
Since $\sigma(1)$ and $\sigma(2)$ are connected and $\sigma(2)>\sigma(1),$ then $\sigma(2)=i_1+1.$ By similar arguments, 
$$\sigma(3)=i_1+2,\ldots,\sigma(k)=i_1+k-1,$$
but $\sigma(k+1)\ne i_1+k.$ Such situation is possible only in one case:  
$i_1=~2n-~k+~1$ and $ \sigma(k+~1)~=~1.$ So,
$$\sigma(k+1)=1, \sigma(k+2)=2,\ldots,\sigma(2n)=i_1-1.$$
In other words,
$$\sigma=i_1\;(i_1+1)\;\ldots\; (2n)\;1\;2\;\ldots (i_1-1).$$
We obtain permutation $\sigma$ that has exactly one up-run if $i_1=1,$ 
or two up-runs if $i_1>1.$ 
So, we obtain permutations of type {\bf I} or {\bf III}.
Therefore, $\sigma\in D_{2n}.$

{\bf Case B.} Now consider the case $\sigma(1)=i_1>\sigma(2).$ Take $k>1,$ such that 
$\sigma(k-1)>\sigma(k)$ and $\sigma(k+1)>\sigma(k).$ 

Since $\sigma(1)$ and $\sigma(2)$ are connected and $\sigma(2)<\sigma(1),$ then $\sigma(2)=i_1-1.$ By similar arguments, 
$$\sigma(3)=i_1-2,\ldots,\sigma(k)=i_1-k+1.$$
but $\sigma(k+1)\ne i_1-k.$ Such situation is possible only in one case:  
$i_1=k, \sigma(k+1)=2n.$ So,
$$\sigma(k+1)=2n, \sigma(k+2)=2n-1,\ldots,\sigma(2n)=i_1+1.$$
In other words,
$$\sigma=i_1\;(i_1-1)\;\ldots\; 1\;{2n}\; (2n-1) \;\ldots (i_1+1)$$
We obtain permutation $\sigma$ that has exactly one down-run if $i_1=2n$ or two down-runs if $i_1<2n.$
In other words we obtain permutations of type {\bf II} or {\bf IV}.
Thus, $\sigma\in D_{2n}.$
$\square$

\section{Proof of Theorem \ref{main}}

\sudda{Main result of our paper is the following 
\begin{thm} \label{main} If $A=(a_{i,j})_{1\le i,j\le 2n}$ is symmetric, then symmetry group of the pfaffian polynomial $pf_{2n}= pf_{2n} \bar A$ is isomorphic to dihedral group,
$$Sym\, pf_{2n}\cong D_{2n}.$$
\end{thm}
}

First we prove that $D_{2n}\subseteq Sym\;pf_{2n}.$

\begin{lm} \label{sympfaff}
If $A=(a_{i,j}) $ is symmetric, then pfaffian is invariant under action of dihedral group $D_{2n},$ 
$$\mu (pf_{2n})= pf_{2n}$$
for any $\mu \in D_{2n}.$
\end{lm}

{\bf Proof.} The dihedral group $D_{2n}$ has order $4n$ and is generated by cyclic permutation 
$$\sigma=\left(\begin{array}{cccccc}
1&2&3&\cdots &2n-1&2n\\
2&3&4&\cdots &2n&1\end{array}\right)$$
and reflection 
$$\tau=
\left(\begin{array}{cccccccccc}
1&2&3&\cdots &n&n+1&n+2&\cdots   &2n-1&2n\\
1&2n&2n-1&\cdots &n+2&n+1&n&\cdots &3&2\end{array}\right)$$

To prove our lemma it is enough to establish that 
$$\sigma (pf_{2n})=pf_{2n},$$
$$\tau ( pf_{2n})=pf_{2n},$$
if $a_{i,j}=a_{j,i},$ for any $1\le i<j\le 2n.$ 

 Recall that 
 $\alpha=(i_1,i_2,\ldots,i_{2n-1},i_{2n})$ is Pfaff permutation, if 
$$i_1<i_3<i_5<\cdots <i_{2n-1},$$
$$i_1<i_2,i_3<i_4,\ldots, i_{2n-1}<i_{2n}.$$
Let $S_{2n,pf}$ be set of Pfaff permutations. 
Below we use one-line notation for permutations. We write 
 $\alpha=(i_1,i_2,\ldots,i_{2n-1},i_{2n})$ 
 instead of 
 $$\alpha=\left(\begin{array}{ccccc}
 1&2&\cdots&2n-1&2n\\
  i_1&i_2&\ldots &i_{2n-1}&i_{2n}
  \end{array}\right).$$ 

Note that 
$$\tau(i)+i=\left\{\begin{array}{ll} 2n+2, & \mbox{ if $1<i\le 2n.$}\\
2& \mbox{ if $i=1$}\end{array}\right.$$
Set
$$\bar i=2n+2-i,$$
if $i>1.$

Now we study action of the generator $\sigma$ on pfaffian polynomials, when generators are symmetric, $a_{i,j}=a_{j,i},$ for any $1\le i,j\le 2n.$  

Let  $\alpha= (1,i_2,i_3,\ldots,i_{2n})\in S_{2n,pf},$ and $l=\alpha^{-1}(2n).$ Then $l$ is even, $l=2k,$ and 
$$\sigma(a_{\alpha})=\sigma(a_{i_1,i_2}\cdots a_{i_{2n-1}i_{2n}})= $$

$$a_{i_1+1,i_2+1}\cdots a_{i_{2k-3}+1,i_{2k-2}+1} a_{i_{2k-1}+1,1}a_{i_{2k+1}+1,i_{2k+2}+1} \cdots a_{i_{2n-1}+1,i_{2n}+1} =
$$

$$a_{\tilde \alpha},$$
where
$$\tilde\alpha=(1,i_{2k-1}+1, i_1+1,i_2+1,
\ldots, {i_{2k-3}+1,i_{2k-2}+1}, {i_{2k+1}+1,i_{2k+2}+1} , \ldots ,{i_{2n-1}+1,i_{2n}+1} ).$$
Here we change $a_{i_{2k-1}+1,1}$ to $a_{1,i_{2k-1}+1}.$
We see that the map
$$S_{2n,pf}\rightarrow S_{2n,pf}, \quad \alpha\mapsto \tilde\alpha$$
is bijection and 
$$sign\; \tilde\alpha= sign\; \alpha.$$

Hence
$$\sigma(pf_{2n})=\sum_{\alpha\in S_{2n,pf}} sign\,\alpha\,\sigma(a_{\alpha})=
 \sum_{\alpha\in S_{2n,pf}} sign\,\tilde\alpha\,a_{\tilde \alpha}=pf_{2n}.$$
So, we have estabilshed that pfaffian is invariant under action $\sigma\in D_{2n}. $

Let us study action of the generator $\tau$ on pfaffian polynomials. 

We have
$$\tau: a_{\alpha}\mapsto a_{1,\;\overline{i_2}}a_{\overline{i_3},\; \overline{i_4}}\cdots a_{\overline{i_{2n-1}},\;\overline{i_{2n}}}.$$
Since
$$\overline{i_{2k-1}}>\overline{i_{2k}}, \quad 1<k\le n,$$
we have to change 
$a_{\overline{i_{2k-1}},\; \overline{i_{2k}}}$ to $a_{\overline{i_{2k}},\; \overline{i_{2k-1}}}.$ 
Further,
$$\overline{i_3}>\overline{i_5}>\cdots >\overline{i_{2n-1}}>1.$$
Therefore, 
$$\tau:a_\alpha\mapsto a_{\bar\alpha},$$
where
$$a_{\bar\alpha}=
a_{1,\; \overline{i_2}} a_{\overline{i_{2n}},\; \overline{i_{2n-1}}}a_{\overline{i_{2n-2}},\; \overline{i_{2n-3}}}\cdots a_{\overline{i_4},\; \overline{i_3}}.$$

We see that 
$$sign\; \alpha= sign \; \bar\alpha.$$
Note that the map
$$S_{2n,pf}\rightarrow S_{2n,pf}, \quad \alpha\mapsto \bar\alpha,$$
is bijection.
Therefore,
$$\tau(pf_{2n})=\sum_{\alpha\in S_{2n,pf}} sign\, \alpha \; \tau(a_{\alpha})= 
\sum_{\alpha\in S_{2n,pf}} sign\,\bar\alpha \; a_{\bar\alpha}= pf_{2n}.$$
So, we have proved that  pfafian $pf_{2n}$ is invariant under action of dihedral group $D_{2n} $ of order $4n,$ if the matrix $(a_{i,j})_{1\le i,j\le 2n}$ is symmetric. 
$\square$

{\bf Example.} Let 
$$\tau=
\left(\begin{array}{cccc}1&2&3&4\\
1&4&3&2\end{array}\right), \mu=\left(\begin{array}{cccc}1&2&3&4\\
4&3&2&1\end{array}\right).$$
Then
$$\tau (pf_4)=\tau(a_{1,2}a_{3,4}-a_{1,3}a_{2,4}+a_{1,4}a_{2,3})=$$
$$a_{1,4}a_{3,2}-a_{1,3}a_{4,2}+a_{1,2}a_{4,3}=$$
$$a_{1,4}a_{2,3}-a_{1,3}a_{2,4}+a_{1,2}a_{3,4}=pf_4,$$

$$\mu (pf_4)=\mu(a_{1,2}a_{3,4}-a_{1,3}a_{2,4}+a_{1,4}a_{2,3})=$$
$$a_{4,3}a_{2,1}-a_{4,2}a_{3,1}+a_{4,1}a_{3,2}=$$
$$a_{3,4}a_{1,2}-a_{2,4}a_{1,3}+a_{1,4}a_{2,3}=pf_4.$$

{\bf Proof of Theorem \ref{main}}. Let $\sigma\in Sym\,pf_{2n}$ i.e.,
$$\sigma\,(pf_{2n})=pf_{2n}$$
for any $a_{i,j},$ such that $a_{i,j}=a_{j,i}.$
In particular, $\sigma$ is a symmetry of the pfaffian polynomial $pf_{2n} ((x_i-x_j)^2)_{1\le i<j\le 2n}.$
By Theorem \ref{translation} the polynomial $g_{2n}(x_1,\ldots,x_{2n})$ is pfaffian polynomial $pf_{2n} ((x_i-x_j)^2)_{1\le i<j\le 2n}$ up to non-zero constant. 
Therefore by  Lemma \ref{sym g} $\sigma\in D_{2n}.$

In remains to use Lemma \ref{sympfaff} to finish the proof. $\square$

\section{ Proof of Theorem \ref{cos}}

 Let $a_{i,j}=\cos(x_i-x_j).$ We have to prove that
$$pf_{2n}=\cos(x_1-x_2+x_3-x_4+\cdots+x_{2n-1}-x_{2n}).$$
The proof is based on the following elementary trigonometric facts.
  
 \begin{lm}\label{17jan2022} 
$$-\cos \alpha \cos (\theta-\alpha)+\cos \beta \cos(\theta-\beta)=\sin(\alpha-\beta)\sin(\alpha+\beta-\theta).$$
\end{lm}

{\bf Proof.} We have 
$$\sin(\alpha-\beta)\sin(\alpha+\beta-\theta)=$$
$$\sin(\alpha-\theta+\theta-\beta)\sin(\alpha+\beta-\theta)=$$

$$(\sin(\alpha-\theta)\cos(\theta-\beta)+
\cos(\alpha-\theta)\sin(\theta-\beta))
\sin(\alpha+\beta-\theta)=$$

{$$\sin(\alpha-\theta)\cos(\theta-\beta)\sin(\alpha+\beta-\theta)+
\cos(\alpha-\theta)\sin(\theta-\beta)
\sin(\alpha+\beta-\theta)=$$}

$$\sin(\alpha-\theta)\cos(\theta-\beta)\sin(\alpha+\beta-\theta)+
\cos(\alpha-\theta)\cos(\theta-\beta)
\cos(\alpha+\beta-\theta)+$$
$$
\cos(\alpha-\theta)\sin(\theta-\beta)
\sin(\alpha+\beta-\theta)-
\cos(\alpha-\theta)\cos(\theta-\beta)
\cos(\alpha+\beta-\theta)
=$$

$$\cos(\theta-\beta)\left\{\sin(\alpha-\theta)\sin(\alpha+\beta-\theta)+
\cos(\alpha-\theta)\cos(\alpha+\beta-\theta)\right\}+$$
$$
\cos(\alpha-\theta)\{\sin(\theta-\beta)
\sin(\alpha+\beta-\theta)-
\cos(\theta-\beta)
\cos(\alpha+\beta-\theta)\}
=$$

$$\cos(\theta-\beta)
\cos(\alpha-\theta-(\alpha+\beta-\theta))$$
$$-
\cos(\alpha-\theta)
\cos(\theta-\beta+(\alpha+\beta-\theta))
=$$

$$\cos(\theta-\beta)
\cos \beta-\cos(\theta-\alpha)\cos(\alpha).$$
$\square$

\begin{lm}\label{29.04.2019}
$$\sum_{i=1}^n (-1)^i \sin\alpha_i \, \sin( \sum_{j=1}^{i-1}(-1)^j \alpha_j-\sum_{j=i+1}^n (-1)^j \alpha_j)=0,$$

$$\sum_{i=1}^n (-1)^i \cos\alpha_i \, \sin( \sum_{j=1}^{i-1}(-1)^j \alpha_j-\sum_{j=i+1}^n (-1)^j \alpha_j)=
\left\{\begin{array}{cc} 0&\mbox{ if $n$ is odd}\\
-\sin (\sum_{j=1}^n (-1)^j \alpha_j)&\mbox{ if $n$ is even}\end{array}\right.$$
\end{lm}

{\bf Proof.}  Follows  from known relation
$$\sin\, \sum_{j=1}^n \theta_i =\sum_{\mbox{\small odd $k\ge 1$}}
 (-1)^{\frac {(k - 1)} {2}}\sum_ {A\subseteq [n], | 
      A | = k}\; (\prod_ {i\in A}\sin\theta_i \prod_ {i \notin A} \cos\theta_i).
$$
$\square$

We use induction on $n.$ For $n=2$ our statement is evident. Suppose that for $n-1$ our statement is true,
$$ pf_{2n-2}(2,\ldots,\hat{i},\ldots,2n)=\cos (x_2-x_3+\cdots +(-1)^{i-1}x_{i-1}+(-1)^i x_{i+1}+\cdots
 + x_{2n-1}-x_{2n}).$$

Let us decompose pfaffian by first row
$$pf_{2n}=\sum_{i=2}^{2n} (-1)^i a_{1,i}\; pf_{2n-2}(2,\ldots,\hat{i},\ldots,2n).$$
Therefore, 
$$pf_{2n}=\sum_{i=2}^{2n} (-1)^i \cos(x_1-x_i)\; \cos (x_2-x_3+\cdots +(-1)^{i-1}x_{i-1}+(-1)^i x_{i+1}+\cdots
 + x_{2n-1}-x_{2n})=$$
$$R_1+R_2,$$
where
$$R_1=$$
$$-\cos(x_1-x_{2n-1})\cos(x_2-x_3+\cdots +x_{2n-2}-x_{2n})+\cos(x_1-  x_{2n})\cos(x_2-x_3+\cdots+x_{2n-2}-x_{2n-1}),$$

$$R_2=$$ $$
 \sum_{i=2}^{2n-2} (-1)^i \cos(x_1-x_i) \; \cos (x_2-x_3+\cdots +(-1)^{i-1}x_{i-1}+
 (-1)^{i}x_{i+1}+\cdots + x_{2n-1}-x_{2n}).$$
   
We have   
   $$R_2=
 \sum_{i=2}^{2n-2} (-1)^i \cos(x_1-x_i) \;\cos(x_{2n-1}-x_{2n})R_{2,i}'-(-1)^i \cos(x_1-x_i)\sin(x_{2n-1}-x_{2n}) R_{2,i}'',$$
 where 
$$R_{2,i}'=
  \cos (x_2-x_3+\cdots +(-1)^{i-1}x_{i-1}+
 (-1)^{i}x_{i+1}+\cdots+x_{2n-3}-x_{2n-2} ),$$ 
 $$
 R_{2,i}''=
  \sin (x_2-x_3+\cdots +(-1)^{i-1}x_{i-1}+
 (-1)^{i}x_{i+1}+\cdots +x_{2n-3}-x_{2n-2}).
 $$

By Lemma \ref{29.04.2019}
$$\sum_{i=2}^{2n-2} 
-(-1)^i \cos(x_1-x_i)\sin(x_{2n-1}-x_{2n}) R_{2,i}''=$$

$$-\sin(x_{2n-1}-x_{2n})\sum_{i=2}^{2n-2} 
(-1)^i \cos(x_1-x_i) 
 \sin (x_2-x_3+\cdots +(-1)^{i-1}x_{i-1}+
 (-1)^{i}x_{i+1}+\cdots +x_{2n-3}-x_{2n-2})
=0,$$
and
      $$R_2=\cos(x_{2n-1}-x_{2n})\;\sum_{i=2}^{2n-2} (-1)^i \cos(x_1-x_i)R_{2,i}'=$$
      
   $$\cos(x_{2n-1}-x_{2n})    \times $$ $$
 \sum_{i=2}^{2n-2} (-1)^i \cos(x_1-x_i) \cos (x_2-x_3+\cdots +(-1)^{i-1}x_{i-1}+
 (-1)^{i}x_{i+1}+\cdots + x_{2n-3}-x_{2n-2}).
  $$
  Therefore by inductive suggestion,
  $$R_2=\cos(x_{2n-1}-x_{2n})\; pf_{2n-2}(1,2,\ldots,2n-2) =\cos(x_{2n-1}-x_{2n})\cos(\sum_{i=1}^{2n-2}(-1)^j x_j).$$

By Lemma \ref{17jan2022} applied  for
 $$\alpha=x_1-x_{2n-1}, \quad \beta=x_1-x_{2n}, \quad \theta=x_1+x_2-x_3+x_4-\cdots-
 x_{2n-3}+x_{2n-2}-
 x_{2n-1}-x_{2n},$$ we obtain
$$R_1=
\sin(-x_{2n-1}+x_{2n})\sin(x_1-x_2+x_3-x_4+\cdots -x_{2n-2}).$$
  Thus,
  $$pf_{2n}=R_1+R_2=  $$
$$     -\sin(x_{2n-1}-x_{2n})\sin(x_1-x_2+x_3-x_4+\cdots -x_{2n-2})+\cos(x_{2n-1}-x_{2n})\cos(\sum_{i=1}^{2n-2}(-1)^j x_j)=$$
     $$ \cos(\sum_{i=1}^{2n} (-1)^i x_i).$$
$\square$   

{\bf Acknowledgment.}

The work was supported by grant AP08855944.GF  Ministry of Education and Science of Kazakhstan Republic.

\end{document}